\documentclass[11pt,reqno]{amsart}%,monochrome
\usepackage{amssymb}
\usepackage[latin1]{inputenc}
\usepackage{amssymb}
\usepackage[all]{xy}
\usepackage{amsfonts}
\usepackage{float}
\usepackage{dsfont}
\usepackage[frenchb,english]{babel}
\usepackage{enumerate}
\usepackage{amsfonts, amsmath, amssymb, stmaryrd, latexsym}
\usepackage{icomma} % la virgule agira comme sÃÂÃÂ©parateur dÃÂÃÂ©cimal en mode mathÃÂÃÂ©matique
\usepackage{pstricks}
\usepackage[all]{xy}
\usepackage{float}
\usepackage{amscd}
\usepackage{geometry}
\usepackage{centernot}
\usepackage[frenchb,english]{babel}

\begin{document}
	\theoremstyle{definition}
	\newtheorem{defen}{Definition}[]
	\newtheorem{rem}{Remark}[]
	\newtheorem{exam}{Example}[]
	
	\theoremstyle{plain}
	\newtheorem{them}{Theorem}[]
	\newtheorem{prop}{Proposition}[]
	\newtheorem{cor}{Corollary}[]
	\newtheorem{lem}{Lemma}[]
	
	%%%%%%%%%%%%%%%%%%%%%%%%%%%%%%%%%%%%%%
	%%%%%%%%%%%%%%%%%%%%%%%%%%%%%%%%%%%
	\def\NN{\mathds{N}}
	\def\RR{\mathbb{R}}
	\def\HH{I\!\! H}
	\def\QQ{\mathbb{Q}}
	\def\CC{\mathbb{C}}
	\def\ZZ{\mathbb{Z}}
	
	\makeatletter
	\@namedef{subjclassname@2020}{\textup{2020} Mathematics Subject Classification}
	\makeatother
	
	\title[P\'{o}lya group and Capitulation...]{P\'{o}lya groups and capitulation in some real biquadratic number fields}
	
	\author[Said El madrari]{S. EL MADRARI}
	\author[Mohammed Taous]{M. Taous}
	\author[Abdelkader Zekhnini]{A. Zekhnini$^{\ast}$}
	
	\address{Said El madrari, Department of  Mathematics, Moulay Ismail University, Faculty of  Sciences and Technology, Errachidia, Morocco.}
	\email{saidelmadrari@gmail.com}
	\address{Mohammed Taous, Department Mathematics, Mohammed First University, Faculty of  Sciences, Oujda, Morocco.}
	\email{taousm@hotmail.com}
	
	\address{Abdelkader Zekhnini, Department Mathematics, Mohammed First University, Faculty of  Sciences, Oujda, Morocco}
	\email{zekha1@yahoo.fr}
\thanks{$^\ast$Corresponding author.}
	\begin{abstract}
	The P\'olya group, $\mathcal{P}_O(K)$, of a  number field $K$ is the subgroup of its ideal class group $\mathrm{C}_K$ generated by the classes of the products of the maximal ideals of $K$ having  the same absolute norm. $K$ is said to be a P\'olya field if $\mathcal{P}_O(K)$ is trivial. In this note, we are interested in  determining, using capitulation theory, the  P\'{o}lya groups of some real biquadratic number fields $K=\mathbb{Q}(\sqrt{d_{1}},\sqrt{d_{2}})$, where  $d_1$ and $d_2$ are two square-free real integers satisfying certain conditions. As a consequence, we   deduce all such  P\'{o}lya fields  $K$,  and the structure of their first cohomology group of units.
	\end{abstract}

	\maketitle
	\section{\bf Introduction}\label{section1}
	Let $K$ be an algebraic number field. Denote by $\mathcal{O}_K$ and $\mathrm{C}_K$  its ring of integers and its ideal class group respectively.  For each integer $q\geq 2$, let  $\prod_q(K)$ be the product of all  maximal ideals of $\mathcal{O}_K$ with the same absolute norm  $q$,  if these ideals do not exist we set  $\prod_q(K)=\mathcal{O}_K$. In their monograph \cite[\S II.4]{CahCha} about the ring of \textit{integer-valued polynomials} of $\mathcal{O}_K$, that is  $\mathrm{Int}(\mathcal{O}_K)=\{P\in K[X]\;|\;P(\mathcal{O}_K)\subset \mathcal{O}_K\}$,  Cahen and  Chabert introduced in 1997 the notion of {\it P\'{o}lya group}, $\mathcal{P}_O(K)$, of $K$  which is the  subgroup of $\mathrm{C}_K$  generated by the classes of the ideals $\prod_q(K)$. This group is then  used to characterize a class of number fields called {\it P\'olya fields} \cite{Ost}, these are fields $K$ such that  $\mathcal{P}_O(K)$ is trivial (equivalently  the ideals $ \prod_{q}(K)$ are principal). Note also that  the notion of P\'{o}lya  fields,  firstly introduced by Zantema in \cite{Za}, follows from the works of P\'{o}lya \cite{Pol} and Ostrowski \cite{Ost} on the integer-valued polynomials. They were interested in whether $\mathrm{Int}(\mathcal{O}_K)$ has a regular basis  as a free $\mathcal{O}_K$-module (see \cite{CahCha}). According to Zantema \cite{Za}, a field $K$ is  a {\it P\'olya field} if $\mathrm{Int}(\mathcal{O}_K)$ has a regular basis $(e_n)_{n\in\NN}$, i.e., the degree of each polynomial $e_n$ is $n$ (see \cite{Pol}).
	
	Quadratic P\'olya fields are well known. They  were completely determined by Zantema in \cite{Za}.  They can also be deduced   from the ambiguous ideals (See \cite[\S75]{Hil95} and \cite[pages: 41-42]{CahCha}). After quadratic fields,  it is natural that many mathematicians are interested in biquadratic number fields that are P\'olya fields. Thereby,   A. Leriche \cite{Leri1, Leri2}  gave some biquadratic  P\'{o}lya fields. In \cite{CM, TauZekh, Ze-16},  the last two authors of the present paper, characterize all imaginary  biquadratic P\'{o}lya   fields by  firstly determining their P\'{o}lya groups. B. Heidaryan,  A. Rajaei  in $2017$ (\cite{HeiRaj17}) and  Tougma in $2021$ (\cite{Tougma})   gave partial  results about  real biquadrtic P\'{o}lya  fields.  In this note and by using capitulation theory, we are interested in characterizing some  real biquadrtic P\'{o}lya fields.
	
Let $d_1$ and $d_2$ be two square-free real coprime integers such that  $d_1$ or $d_2\equiv1\pmod4$, and $d_3=d_1d_2$ is divisible by a prime congruent to $3 \pmod4$ or the norm of the fundamental unit of $k=\mathbb{Q}(\sqrt{d_{3}})$ equals  $-1$. Put $K=\QQ(\sqrt{d_1},\sqrt{d_2})$, in this paper, we aim  to characterize  the P\'{o}lya group of the real biquadratic number field $K$. Precisely we evaluate the P\'{o}lya group of $K$ in terms of that of $k$ and of the  classes of $k$  capitulating in $K$. Next, we deduce all real biquadratic P\'{o}lya fields  $K$ satisfying these conditions, and thus Theorems $A$ and $D$ of  \cite{HeiRaj17} turn out to be  special cases of Theorem \ref{19} (see below)  which deals with general biquadratic number fields (for more details, see remark \ref{23} below). As a consequence, we deduce the structure of the first cohomology group of units of  $K$.

	\section*{\bf Notations}
	Let $d_1>1$ and $d_2>1$ be two square-free coprime integers, put $d_3=d_1d_2$. The following notations will be used through this paper:
	\begin{itemize}
		\item $K=\QQ(\sqrt{d_1},\sqrt{d_2})$:  a real biquadratic number field.
		\item $e_p$:   the ramification index of  a prime number $p$ in $K/\QQ$.
		\item $\epsilon_i$:  the fundamental unit of  $k_i=\QQ(\sqrt{d_i})$, $i=1, 2, 3$, some times we use $k$ as $k_3$.
		\item $N(\epsilon_i)=N_{k_i/\mathbb{Q}}(\epsilon_i)$ : the absolute norm  of $\epsilon_i$.
		\item $N_{L/F}$ : the relative norm of the extension $L/F$.
		\item $L^*$: the multiplicative group  of nonzero elements of a number field $L$.
		\item $\mathcal{O}_L$:  the ring of integers of a number field $L$.
		\item $E_L$:  the unit group of a number field $L$.
		\item $G_L$:  the Galois group of a number field $L$ over $\QQ$.
		\item  $H^1(G_L, E_L)$:   the first cohomology group of units of $L$.
		\item $d_L$:   the discriminant of a number field $L$.
		\item $Q=\left[E_K: \prod_{i=1}^3 E_{k_i}\right]$ : the unit index of $K$.
		\item[$\bullet$] $\left(\frac{m }{n}\right)$: the Legendre symbol.        			
		%\item[$\bullet$] $\left(\frac{m }{n}\right)_{4}$ : the rational biquadratic symbol. If $n\equiv 1\pmod 8$,  $\left(\frac{n }{2}\right)_{4}:=(-1)^{\frac{n-1}{8}}$.
		\item $E_b=(\mathbb{Z}/2\mathbb{Z})^b$, with $b\in\NN^*$.	
	\end{itemize}
	\section{\bf Preliminaries}
	Let $F$ be a number field and $L/F$ a quadratic extension. Let $Am(L/F)$ denote the ambiguous class group of $L/F$, that is the group of classes  $c\in C_L$ fixed by a generator $\sigma$ of the Galois group of $L/F$. If the class number of $F$ is odd, then $Am(L/F)$ is exactly  the group of classes $c$ of order $\leq 2$ (see \cite[Proposition 2.4]{Mo16}).  The classes of ideals which are invariant by $\sigma$  constitute a subgroup of $Am(L/F)$ called  the group of strongly ambiguous classes of $L/F$, and it will be denoted $Am_s(L/F)$. The following result gives the order of   $Am(L/F)$ and $Am_s(L/F)$ (e.g., \cite{Lemm-13}).
	\begin{them}\label{1}
		Let $L/F$ be a  quadratic extension. Then the number of ambiguous classes is given by	
		\begin{equation*}
		\begin{split}
		|Am(L/F)|&=h(F)\frac{2^{t-1}}{[E_{F}:E_{F}\cap N_{L/F}(L^{*})]}\\
		&=|Am_s(L/F)|\times[E_{F}\cap N_{L/F}(L^{*}) :E_{F}\cap N_{L/F}(E_L)],
		\end{split}
		\end{equation*}
		where $t$ is the number of finite and infinite primes of $F$ ramified in $L$, and $h(F)$ the class number of the field $F$.
	\end{them}
	\begin{cor}\label{cor1}
		Let $L=\QQ(\sqrt{d})$ be a quadratic number field, where $d$ is a square-free integer,  and $\epsilon$ be its fundamental unit. Then $Am(L/\QQ)=Am_s(L/\QQ)$ if and only if $N(\epsilon)=-1$  or there exists a prime integer $p \mid d$ such that $p\equiv 3\pmod4$.
		
	\end{cor}
	\begin{proof}
		As $L/\QQ$ is a quadratic extension,  so by  Theorem \ref{1} we get:
		$$|Am(L/\QQ)|=h(\QQ)\frac{2^{t-1}}{[E_{\QQ}:E_{\QQ}\bigcap N_{L/\QQ}(L^{*})]}\cdot$$
		Note that the infinite prime of $\QQ$ does not ramify in $L$, so $t$ is the number of prime integers that ramify in $L$.
		As $E_{\QQ}=\langle -1\rangle$, so to compute the index $[E_{\QQ}:E_{\QQ}\bigcap N_{L/\QQ}(L^{*})]$ we must study in which cases there exists $x\in L^*$ such that  $N_{L/\QQ}(x)=-1$. Note also that  $h(\QQ)=1$ and
		$[E_{\QQ}:N_{L/\QQ}(E_L)]=[E_{\QQ}:E_{\QQ}\cap N_{L/\QQ}(L^{*})][E_{\QQ}\cap N_{L/\QQ}(L^{*}):N_{L/\QQ}(E_L)]\cdot$
		Thus we will discuss the following cases.
		\begin{enumerate}[$\bullet$]
			\item If $N_{L/\QQ}(\epsilon)=-1$,  then $-1\in N_{L/\QQ}(E_L)$. As  $E_L=\langle -1,  \epsilon\rangle$, so  $[E_{\QQ}:N_{L/\QQ}(E_L)]=1$,  then $[E_{\QQ}:E_{\QQ}\cap N_{L/\QQ}(L^{*})]=[E_{\QQ}\cap N_{L/\QQ}(L^{*}):N_{L/\QQ}(E_L)]=1$,  by Theorem \ref{1} $|Am(L/\QQ)|=|Am_s(L/\QQ)|=2^{t-1}$. Then $Am(L/\QQ)=Am_s(L/\QQ)$.
			\item If there exists a prime $p |d$ such that $p\equiv 3\pmod4$,  then $N_{L/\QQ}(\epsilon)=1$,  so $[E_{\QQ}:E_{\QQ}\cap N_{L/\QQ}(L^{*})]=2$ and $[E_{\QQ}\cap N_{L/\QQ}(L^{*}):N_{L/\QQ}(E_L)]=1$,  consequently $|Am(L/\QQ)|=|Am_s(L/\QQ)|=2^{t-2}$.
			Then $Am(L/\QQ)=Am_s(L/\QQ)$.
			\item Assume $N_{L/\QQ}(\epsilon)=1$  and	any odd prime $p$ dividing $d$ satisfies $p\equiv 1\pmod4$. Denote, for all $p|d$, by $\mathfrak{p}$ a prime ideal of $L$ above $p$, then Hilbert symbol implies that $\left(\dfrac{-1, d_3}{\mathfrak{p}}\right)=\left(\dfrac{-1}{p}\right)=1$. Hence the equation $N_{L/\QQ}(x)=-1$ has a solution in $L^*$, which implies that $[E_{\QQ}\cap N_{L/\QQ}(L^{*}):N_{L/\QQ}(E_L)]=2$. we then deduce that $\left|Am(L/\QQ)\right|=2\left|Am_s(L/\QQ)\right|$.
			
		\end{enumerate}
	\end{proof}
	
	The first important result about P\'olya groups was given by Hilbert \cite[\S75]{Hil95} who studied the ambiguous ideals, since, for a quadratic number field $L$, the P\'olya group of $L$ is  the group generated by the strongly ambiguous classes, thus $\mathcal{P}_O(L)=Am_s(L/\QQ)$.
	\begin{prop}\label{10}
		Let $L$ be a quadratic number field and $\epsilon$ its fundamental unit $($if L is real$)$. If $t$ denotes the number of ramified primes  in the extension $L/\QQ$, then
		$$\left|\mathcal{P}_O(L)\right|=
		\left\{\begin{array}{ll}
			2^{t-2} & \text{ if $L$ is real and } N_{L/\QQ}(\epsilon)=1,\\
			2^{t-1} & \text{ otherwise}.
		\end{array}\right.
		$$
	\end{prop}
	\noindent Hence one can   easily deduce  all quadratic P\'olya number fields.
	\begin{prop}[\cite{Za}]
		A quadratic field $L=\QQ(\sqrt d)$, for a square-free integer $d$, is a P\'olya field if and
		only if $d$ is  one of the following forms, where
		$p$ and $q$ denote two distinct odd prime numbers:
		\begin{enumerate}[\rm1.]
			\item $d = -2$, or $-1$, or $2$, or $ -p$ with $p \equiv 3 \pmod 4$, or  $ p$.
			\item $d = 2p$ and
			$\begin{array}[t]{ll}
				\text{-  either } p \equiv 3 \pmod 4\\
				\text{- or } p \equiv 1 \pmod 4 \text{ and the norm of the fundamental unit of $L$ is $1$}.
			\end{array}$
			\item $d = pq$ and
			$\begin{array}[t]{ll}
				\text{- either } p, q \equiv 3 \pmod 4 \\
				\text{- or } p, q \equiv 1 \pmod 4 \text{ and  the norm of the fundamental unit of $L$ is $1$}.
			\end{array}$
		\end{enumerate}
	\end{prop}
    	Another characterization of P\'{o}lya group of a number field $L$, when $L/\QQ$ is Galois, is thanks to Zantema see \cite[p. 163]{Za}. It gives a cohomological approach of to a P\'{o}lya group.
    \begin{prop}[\cite{Za}]\label{zatm6}
    	Let $L/\QQ$ be a Galois extension and $d_L$ be its
    	discriminant. Denote by $e_p$  the ramification index of  a prime number $p$ in $L$. Then the following sequence is exact $$
    	1\rightarrow H^1(G_L,
    	E_L)\rightarrow\displaystyle\oplus_{p/d_L}\ZZ/e_p\ZZ\rightarrow
    	\mathcal{P}_O(L)\rightarrow 1.
    	$$
    	In particular, $|H^1(G_L, E_L)|\times|\mathcal{P}_O(L)|=\prod_{p|d_L}e_p$.
    \end{prop}
	Let $L/F$ be  a finite extension of number fields, and $j_{F\rightarrow L}:  \mathrm{C}_F\rightarrow \mathrm{C}_L$  the homomorphism induced by extension of ideals from $F$
	to $L$. An  ideal class of $F$ is said to {\it capitulate} in $L$ if it is in $\ker j_{F\rightarrow L}$, the kernel of the	homomorphism $j_{F\rightarrow L}$, which  is called the capitulation kernel.
		
	\begin{prop}\label{zatm9}
		Let $L$ and $F$ be two Galois number fields such that $F \subset L$. If $L/F$ is unramified at all finite primes, then $j_{F\rightarrow L}(\mathcal{P}_O(F))=\mathcal{P}_O(L)$.
	\end{prop}
\begin{proof}
In proposition 3.4 of \cite{Chab03}, Chabert showed that $j_{F\rightarrow L}(\mathcal{P}_O(F))\subseteq\mathcal{P}_O(L)$, then it is easy to see in the proof of this proposition that as  $L/F$ is Galois and unramified at all finite primes, we find
that $j_{F\rightarrow L}(\mathcal{P}_O(F))=\mathcal{P}_O(L)$.
\end{proof}	
	\section{\bf The structure of P\'{o}lya groups of some real bicyclic biquadratic  number fields}
	In this section, we will characterize some P\'{o}lya groups of the real biquadratic number fields $K=\QQ(\sqrt{d_1},\sqrt{d_2})$, assuming  $K$ is an unramified quadratic extension of $k=\QQ(\sqrt{d_3})$ (see Remark \ref{2} below). All the following results are corollaries  of the  following theorem:
\subsection{Main result}
	\begin{them}\label{3}
	Let $K=\QQ(\sqrt{d}_{1}, \sqrt{d}_{2})$, where $d_1$ and $d_2$  are two  coprime square-free integers. Put $k=\QQ(\sqrt{d_3})$, with $d_3=d_1d_2$, and  assume  $d_1$ or $d_2\equiv1\pmod4$. Put  $j:=j_{k\rightarrow K}$. Let $\epsilon_3$ be the fundamental unit of $k$. Then $\mathcal{P}_O(K)\simeq \mathcal{P}_O(k)/\ker(j)$ if one of  the following cases holds:
	\begin{enumerate}[\rm a.]
		\item $N(\epsilon_{3})=-1$,
		\item $d_3$ is divisible by a prime congruent to $3\pmod4$.
	\end{enumerate}
\end{them}
\begin{rem}\label{2}
	The conditions of the Theorem \ref{3} are equivalent to $K/k$ is  an unramified quadratic extension (\cite{BlaSpeKen}).
\end{rem}
\begin{rem}
	This theorem transforms the problem of determining $\mathcal{P}_O(K)$ to a problem of capitulation, since  the P\'olya group of a quadratic field is well known  (see Proposition \ref{10}).
\end{rem}
\begin{proof}[\bf Proof of Theorem \ref{3}]
	Using some results from the previous section, we are able to prove  this theorem. In fact, this result is the key of all the results of this paper and establishes a relationship between  $\mathcal{P}_O(K)$, $\mathcal{P}_O(k)$ and the capitulation kernel $\ker(j)$. As the reader will remark, this theorem transforms the problem of determining $\mathcal{P}_O(K)$ to a problem of capitulation because  the P\'olya group of a quadratic field is known since Hilbert (see Proposition \ref{10}).
	To continue, recall first that, as $K/k$ is unramified, so   Remark \ref{2} and   Proposition \ref{zatm9} give that  $j(\mathcal{P}_O(k))=\mathcal{P}_O(K)$. On the other hand, it is known by definition that $\mathcal{P}_O(k)\subseteq \mathrm{C}_k$,  and $\mathcal{P}_O(K)\subseteq \mathrm{C}_K$. Let $\theta=j_{\mathcal{P}_O(k)}$ denote the restriction of $j$ to $\mathcal{P}_O(k)$,  then
	$\theta: \mathcal{P}_O(k)\longrightarrow \mathcal{P}_O(K)$ is a group homomorphism,  so $\mathcal{P}_O(K)=j(\mathcal{P}_O(k))\cong \mathcal{P}_O(k)/\ker\theta$,  where of course $\ker\theta=\ker(j)\cap \mathcal{P}_O(k).$\\
	It is well known that $\ker(j) \subseteq Am(k/\QQ)$; indeed as $K/k$ is a quadratic extension, then $N_{K/k}\circ j(c)=c^2$ for all classes $c$ of $\mathrm{C}_k$, so  we deduce that
	 $c  \in \ker(j) $ implies that $c^2=1$, i.e. $c\in Am(k/\QQ)$.\\
	 Note that $\mathcal{P}_O(k) =Am_s(k/\QQ)\subseteq Am(k/\QQ)$. When $N(\epsilon)=-1$  or when there is a prime 	$  p \mid d$ such that $p\equiv 3\pmod4$, Corollary \ref{cor1} yields that  $Am(k/\QQ)=Am_s(k/\QQ)$. From which we deduce that  $\ker(j) \subseteq Am(k/\QQ)=Am_s(k/\QQ)=\mathcal{P}_O(k)$, so  $\ker\theta=\ker(j)\cap \mathcal{P}_O(k)=\ker(j)$. Therefore,  $\mathcal{P}_O(K)\cong \mathcal{P}_O(k)/\ker(j)$. Which completes the proof of the theorem.
\end{proof}
\subsection{\bf P\'{o}lya groups structure}
Keep the notations given in the beginning of  this section and based on a result in \cite{BenSanSny} precisely Theorem \ref{4} below, we will give the structure of P\'{o}lya groups of some fields $K$.
We start by the following theorem which gives the order of $\ker(j)$.
	\begin{them}[\cite{BenSanSny}]\label{4}
		If $K$ is an unramified quadratic extension of a real quadratic number field $k$. Then
		\begin{enumerate}[\rm1.]
			\item $|\ker(j)|=2\Leftrightarrow\left\{\begin{array}{l}
				N(\epsilon_{i})=-1,\text{ for }i=1, 2, 3,\text{ and }\sqrt{\epsilon_{1}\epsilon_{2}\epsilon_{3}}\in K\text{ or }\\
				(N \epsilon_{1}=-1\text{ or } N(\epsilon_{2})=-1), N(\epsilon_{3})=1\text{ and } (\sqrt{\epsilon_{3}}\text{ or } \sqrt{\epsilon_{1}\epsilon_{3}}\text{ or } \sqrt{\epsilon_{2}\epsilon_{3}}\in K).\end{array}\right.$
			\item $|\ker(j)|=8\Leftrightarrow \left\{\begin{array}{ll}
				N(\epsilon_{1})=N(\epsilon_{2})=1\text{ and }N(\epsilon_{3})=-1\text{ or }\\
				N(\epsilon_{i})=1\text{ for } i=1, 2, 3, \text{ and } \sqrt{\epsilon_{1}^{a}\epsilon_{2}^{b}\epsilon_{3}}\notin K\text{ for all }\\
			 a, b\in \{0, 1\}.
			\end{array}\right.$
			\item $|\ker(j)|=4\Leftrightarrow$ anything else occurs.
		\end{enumerate}
	\end{them}
Recall that for  $b\in\NN^*$, $E_b=(\mathbb{Z}/2\mathbb{Z})^b$. We have to discuss two cases. \\
{\bf A. First case: $N(\epsilon_3)=-1$}.
	\begin{them}\label{11}
		Let $K=\QQ(\sqrt{d}_{1}, \sqrt{d}_{2})$, where $d_1$ and $d_2$ are two coprime square-free integers. Put $k=\QQ(\sqrt{d_3})$, where $d_3=d_1d_2$, and assume $N(\epsilon_{3})=-1$. Let  $t$ be  the number of primes  dividing the discriminant $d_{k}$ of $k$.
		Then  $$\mathcal{P}_O(K) \simeq \left\{ \begin{array}{ll}
			E_{t-2}\ \ if \ \ N(\epsilon_{1})=N(\epsilon_{2})=-1\ and\  \sqrt{\epsilon_1\epsilon_2\epsilon_3}\in K, \\
			E_{t-3} \ \  if \ \ N(\epsilon_{1})=N(\epsilon_{2})=-1\ and\  \sqrt{\epsilon_1\epsilon_2\epsilon_3}\notin K\ or\ N(\epsilon_{1})\neq N(\epsilon_{2}), \\
			E_{t-4} \ \ if \ \  N(\epsilon_{1})=N(\epsilon_{2})=1.
			
		\end{array}
		\right.$$	
	\end{them}
	\begin{proof}
		As  $N(\epsilon_{3})=-1$ and $(d_{1}, d_{2})=1$, so $d_1$ or $d_2\equiv1\pmod4$.  Then, by Remark \ref{2}, $K$ is an unramified quadratic extension of $k$. Hence, we are under the conditions of the  Theorem \ref{3}. Since $N(\epsilon_{3})=-1$, so  Proposition \ref{10} implies that  $\mathcal{P}_O(k)\simeq E_{t-1}$.
		We have to distinguish the following cases:
		\begin{enumerate}[\rm1.]
			\item If $N(\epsilon_{1})=N(\epsilon_{2})=-1$ and  $\sqrt{\epsilon_1\epsilon_2\epsilon_3}\in K$, then  by Theorem \ref{4} we get  $|\ker j|=2$, which implies that  $\ker(j)\simeq E_1$. Thus  Theorem \ref{3} yields that  $\mathcal{P}_O(K)\simeq E_{t-1-1}=E_{t-2}$.
			\item If $N(\epsilon_{1})=N(\epsilon_{2})=-1$ and  $\sqrt{\epsilon_1\epsilon_2\epsilon_3}\notin K$, then by Theorem \ref{4} again we get that $|\ker j|=4$; which implies that  $\ker(j)\simeq E_2$. Thus  the Theorem \ref{3} yields that $\mathcal{P}_O(K)\simeq E_{t-1-2}=E_{t-3}$.
			\item If $N(\epsilon_{1})\neq N(\epsilon_{2})$, then by Theorem \ref{4} we get $|\ker(j)|=4$,
			then $\ker(j)\simeq E_2$. Hence $\mathcal{P}_O(K)\simeq E_{t-1-2}=E_{t-3}$.
			\item  If $N(\epsilon_{1})=N(\epsilon_{2})=1$, then by Theorem \ref{4}  we get $|\ker(j)|=8$, and thus $\ker(j)\simeq E_3$. Hence $\mathcal{P}_O(K)\simeq E_{t-1-3}=E_{t-4}$.	
		\end{enumerate}
		Which completes the proof.
	\end{proof}
\noindent {\bf B. Second case:  $d_3$ is divisible by a prime $p\equiv 3\pmod4$}.
	\begin{them}\label{19}
		Let $K=\QQ(\sqrt{d}_{1}, \sqrt{d}_{2})$,  where $d_1$ and $d_2$ are two coprime square-free integers satisfying  $d_1$ or $d_2\equiv1\pmod4$. Let $k=\QQ(\sqrt{d_3})$, with $d_3=d_1d_2$, and assume that $d_3$ is divisible by a prime  $p\equiv 3\pmod4$. Let  $t$ be  the number of primes  dividing the discriminent $d_{k}$ of $k$. Then
		$$\mathcal{P}_O(K)  \simeq \left\{\begin{array}{ll}
			E_{t-3}& if \,  N(\epsilon_{1})\neq N(\epsilon_{2})\ and\  \sqrt{\epsilon_i^a\epsilon_3}\in K \text{ with }a=0, 1; i=1, 2 \text{ such that }  N(\epsilon_{i})=1,   \\
			E_{t-4}&   if \, N(\epsilon_{1})=N(\epsilon_{2})=1\ and\ \sqrt{\epsilon^a_1\epsilon^b_2\epsilon_3}\in K, where \ a,  b \in \{0, 1\},    \\
			& or\ N(\epsilon_{1})\neq N(\epsilon_{2}) \  and\  \sqrt{\epsilon_i^a\epsilon_3}\notin K \text{ with }a=0, 1; i=1, 2 \text{ such that }  N(\epsilon_{i})=1,    \\
			E_{t-5}& if \ N(\epsilon_{1})=N(\epsilon_{2})=1\ and\  \sqrt{\epsilon^a_1\epsilon^b_2\epsilon_3}\notin K,  \ where\ a,  b \in \{0, 1\}.
		\end{array}
		\right.$$
	\end{them}
	\begin{proof}
		As  $d_3$ is divisible by a prime   $p\equiv 3\pmod4$, so $N(\epsilon_{3})=1$.  Since $d_1$ or $d_2\equiv1\pmod4$, then by Remark \ref{2},  $K$ is an unramified quadratic extension of $k$,  hence we are under the conditions of  the Theorem \ref{3}. Moreover, since $N(\epsilon_{3})=1$ and $k$ is real, so  Proposition \ref{10} implies that  $\mathcal{P}_O(k)\simeq E_{t-2}$.
		
		As $p|d_3$,  $(d_1, d_2)=1$ and $d_3=d_1d_2$, so $p$ divides one of the two integers $d_1$ and  $d_2$. Here also we have several cases  to consider:
		\begin{enumerate}[\rm1.]
			\item If $N(\epsilon_{1})=N(\epsilon_{2})=1$ and  $\sqrt{\epsilon^a_1\epsilon^b_2\epsilon_3}\in K$, for some $a,  b\in\{0, 1\}$, then  by Theorem \ref{4} we get $|\ker(j)|=4$,  so $\ker(j)\simeq E_2$. Hence Theorem \ref{3} yields that  $\mathcal{P}_O(K)\simeq E_{t-2-2}=E_{t-4}$.
			\item If $N(\epsilon_{1})=N(\epsilon_{2})=1$ and  $\sqrt{\epsilon^a_1\epsilon^b_2\epsilon_3}\notin K$, for all $a,  b\in\{0, 1\}$, then by Theorem \ref{4} we get $|\ker(j)|=8$,   so $\ker(j)\simeq E_3$. Thus Theorem \ref{3} yields that $\mathcal{P}_O(K)\simeq E_{t-2-3}=E_{t-5}$.
			\item Assume $N(\epsilon_{1})\neq N(\epsilon_{2})$, denote by $\epsilon_i$, with $i=1, 2$,  the fundamental unit such that $N(\epsilon_{i})=1$.
			\begin{enumerate}[\rm a.]
				\item If  $\sqrt{\epsilon_i^a\epsilon_3}\in K$ for some $a=0, 1$,  then by Theorem \ref{4} we get that $|\ker(j)|=2$,  so we have  $\ker(j)\simeq E_1$, and thus $\mathcal{P}_O(K)\simeq E_{t-2-1}=E_{t-3}$.
				\item If  $\sqrt{\epsilon_i^a\epsilon_3}\notin K$ for any  $a=0, 1$,  then by Theorem \ref{4} we get that $|\ker(j)|=4$,  so we have  $\ker(j)\simeq E_2$, and thus $\mathcal{P}_O(K)\simeq E_{t-2-2}=E_{t-4}$.
			\end{enumerate}	
		\end{enumerate}
Which completes the proof.
	\end{proof}
	\begin{rem}\label{23}
		Theorems $A$ and $D$ of  \cite{HeiRaj17} turn out to be  special cases of Theorem \ref{19}  which deals with general cases of biquadratic number fields. Let us clarify this.
		
		Let $K=\QQ(\sqrt{p}, \sqrt{qr})$,  $q\equiv r\equiv-p\equiv1\pmod4$,  $k_1=\QQ(\sqrt{p})$ and $k_2=\QQ(\sqrt{qr})$. Let $\epsilon_{3}=x+y\sqrt{pqr}$ be the fundamental unit of $k=\QQ(\sqrt{d_3})$, where $d_3=pqr$,  and $\epsilon_{1}$, $\epsilon_{2}$ the fundamental units of $k_1$  and $k_2$ respectively. Note that  $N(\epsilon_{3})=1$ and $N(\epsilon_{2})=\pm1$. Since $d_3$ is congruent to $3\pmod4$, then $t=4$. Moreover, $d_3$  is divisible by a prime $p\equiv3\pmod4$, so we can apply Theorem \ref{19}; remark that  Proposition \ref{10} implies that $\mathcal{P}_O(k)\simeq E_{4-2}=E_2$.
		\begin{enumerate}[\rm i.]
			\item Theorem $A$ of  \cite{HeiRaj17} states that if we further assume  $r \equiv 1 \pmod8$, $(\frac{q}{r})=-1$ and $(\frac{p}{r})=1$,  then $\mathcal{P}_O(K) = \ZZ/2\ZZ$. Note first that under these conditions  $N(\epsilon_{1})=1\neq N(\epsilon_{2})=-1$ and $N(\epsilon_{3})=1$. On the other hand, the conditions on $p$, $q$,  $r$ and the unique factorisation in $\ZZ$ imply   that the possible cases of factorization of $x^2-1$ are as follows :
			\begin{center}$\left\{\begin{array}{rl}
					x\pm1 &=y_1^2\\
					x\mp1 &=pqr y_2^2,
				\end{array}\right.$
				$\left\{\begin{array}{rl}
					x\pm1 &=py_1^2\\
					x\mp1 &=qry_2^2,
				\end{array}\right.$  or
				$\left\{\begin{array}{rl}
					x+1 &=2py_1^2\\
					x-1 &=2qry_2^2,
				\end{array}\right.$
			\end{center}
			where $y=y_1y_2$ or $2y_1y_2$. Which implies that $x\pm1$,  $p(x\pm1)$ or $2p(x+1)$ is a square in $\NN$.  Thus according to \cite[Proposition 9]{ATZAR}, $\sqrt{\epsilon_1^a\epsilon_3}\in K$; hence, using Theorem \ref{19}, we deduce that  $\mathcal{P}_O(K)\simeq E_{t-3}=E_1=\ZZ/2\ZZ$  as stated by  Theorem $A$ of  \cite{HeiRaj17}.
			
			More generally, $\mathcal{P}_O(K)$ is always isomorphic to $\ZZ/2\ZZ$, i.e., $K$ is not  a P\'olya field,  whenever  $N(\epsilon_{2})=-1$ and $x\pm1$, $p(x\pm1)$ or $2p(x\pm1)$ is a square in $\NN$ since in these cases $\sqrt{\epsilon_1^a\epsilon_3}\in K$ (see \cite[Proposition 9]{ATZAR}), and  the result comes by Theorem \ref{19}.
			\item Theorem $D$ of  \cite{HeiRaj17} states that  if we further assume  $r \equiv 5 \pmod8$, $(\frac{p}{r})=1$ and $(\frac{p}{q})=-1$,  then $\mathcal{P}_O(K)$ is trivial, i.e., $K$ is a P\'olya field.
			
			Recall   that $N(\epsilon_{1})=N(\epsilon_{3})=1$ and $N(\epsilon_{2})=\pm1$. For $\epsilon_3=x+y\sqrt{pqr}$, and under  the conditions of $p$, $q$,  $r$ we get that
			$x\pm1$,  $p(x\pm1)$ and $2p(x+1)$ are not  squares in $\NN$. Thus according to \cite[Proposition 9]{ATZAR}, if $N(\epsilon_{2})=-1$ one deduces that $\sqrt{\epsilon_1^a\epsilon_3}\notin K$, and according to \cite[Proposition 8]{ATZAR}, if $N(\epsilon_{2})=1$ one deduces that $\sqrt{\epsilon_1^a\epsilon_2^b\epsilon_3}\in K$, which implies, using Theorem \ref{19}, that $\mathcal{P}_O(K)\simeq E_{t-4}=E_0=\{1\}$ as stated by  Theorem $D$ of  \cite{HeiRaj17}.
			
			More generally, $\mathcal{P}_O(K)$ is  always trivial, i.e. $K$ is a P\'olya field,  whenever:
			\begin{enumerate}[$\bullet$]
				\item $N(\epsilon_{2})=1$. In deed, by \cite[Proposition 8]{ATZAR}  we have $\sqrt{\epsilon_1^a\epsilon_2^b\epsilon_3}\in K$. Thus  Theorem \ref{19} implies the result.
				\item $N(\epsilon_{2})=-1$ and $x\pm1$, $p(x\pm1)$ and $2p(x\pm1)$ are not squares in $\NN$ since in this case $\sqrt{\epsilon_1^a\epsilon_3}\notin K$ (see \cite[Proposition 9]{ATZAR}), and  the result is derived by Theorem \ref{19}.
			\end{enumerate}
		\end{enumerate}
	\end{rem}

	\section{\textbf{P\'{o}lya  fields of the form $\QQ(\sqrt{d_1}, \sqrt{d_2})$, where  $(d_1, d_2)=1$ and  $d_1$ or $d_2\equiv1\pmod4$}}
	In this section we will determine all real biquadratic P\'{o}lya fields $K=\QQ(\sqrt{d_1}, \sqrt{d_2})$ such that $(d_1, d_2)=1$ and  $d_1$ or $d_2\equiv1\pmod4$. Put $d_3=d_1d_2$, and let $\epsilon_i$ be the fundamental unit of the quadratic field $\mathbb{Q}(\sqrt{d_i})$, $i=1, 2, 3$. Let $p$, $p_i$ (resp. $q$, $q_i$), where $i\in\NN^*$,  denote odd prime integers congruent to $1\pmod4$ (resp. $3\pmod4$). The two first results in this section are  applications of Theorem \ref{11}. To avoid long theorems, we will consider two cases, in the first one we treat the case $N(\epsilon_{1})= N(\epsilon_{2})=N(\epsilon_{3})=-1$ and in the second one we treat the case $N(\epsilon_{3})=-1$ and ($N(\epsilon_{1})= N(\epsilon_{2})=1$ or $N(\epsilon_{1})\neq N(\epsilon_{2})$).
	\begin{them}
		Let $K=\QQ(\sqrt{d}_{1}, \sqrt{d}_{2})$, where $d_1$ and $d_2$ are  two coprime square-free integers such that $d_1$ or $d_2\equiv1\pmod4$. Let $d_3=d_1d_2$, $Q$  the unit index of $K$ and assume $N(\epsilon_{1})= N(\epsilon_{2})=N(\epsilon_{3})=-1$, let $i\neq j=1, 2$. Then $K$ is a P\'{o}lya field if and only if  one of the following conditions holds.
		\begin{enumerate}[\rm1.]
			\item\label{12}  $d_i=p_1$,  $d_j=p_2$,
			\item $d_i=p$,  $d_j=2$.
			
			\item  $Q=1$  and $d_i$, $d_j$ take the following forms :
			\begin{enumerate}[\rm i.]
				\item\label{13} $d_i=p_1$,  $d_j=p_2p_3$,
				\item\label{14} $d_i=2$, $d_j=p_1p_2$,
				\item\label{15} $d_i=2p_1$, $d_j=p_2$.
			\end{enumerate}
		\end{enumerate}
	\end{them}
	\begin{proof} As $N(\epsilon_{1})= N(\epsilon_{2})=N(\epsilon_{3})=-1$, then according to \cite{Kur-43}, $\{
		\epsilon_1,\epsilon_2, \sqrt{\epsilon_1 \epsilon_2\epsilon_3}\}$ or $\{
		\epsilon_1,\epsilon_2, \epsilon_3\}$ is the fundamental system of units of $K$, it is equivalent to say that $Q=2$ or $1$, depending on whether $\sqrt{\epsilon_1\epsilon_2\epsilon_3}$ is in $K$ or not.  So by Theorem \ref{11} we have $\mathcal{P}_O(K)\simeq E_{t+t_0-3}$, where  $t$ is the number of prime divisors of $d_k$ and $Q=2^{t_0}$. As $N(\epsilon_{3})=-1$, so $d_3$ is not divisible by  any prime congruent to $3\pmod4$. Hence the  prime divisors of the discriminant $d_k$ of $k$ are exactly the prime divisors of  $d_3$. Then
		we have two cases:
		\begin{enumerate}[\rm1.]
			\item If $Q=2$, then $t_0=1$ and $\mathcal{P}_O(K)\simeq E_{t-2}$.  So $K$ is a  P\'{o}lya field if and only if $t=2$. Thus $d_3=p_1p_2$ or $d_3=2p$, i.e., $(d_i=p_1$ and $d_j=p_2)$ or $(d_i=p$ and $d_j=2)$. But with the obtained forms, $Q=2$ is always verified, see \cite[p. 102]{Sc-34}.
			\item If $Q=1$, then $t_0=0$ and $\mathcal{P}_O(K)\simeq E_{t-3}$. So $K$ is a  P\'{o}lya field if and only if $t=3$. Thus $d_3=p_1p_2p_3$ or $d_3=2p_1p_2$.
			Suppose $(d_{i}, d_{j})\equiv(1, 1)\pmod4$,  then  $K$ is a  P\'{o}lya field if and only if  $d_i=p_1$ and $d_j=p_2p_3$.
			Suppose now that $(d_{i}, d_{j})\equiv (2, 1) \pmod4$. Thus $K$ is a  P\'{o}lya field if and only if $(d_i=2$ and $d_j=p_1p_2)$ or $(d_i=2p_1$ and $d_j=p_2)$.
		\end{enumerate}
	\end{proof}

	\begin{them}
		Let $K=\QQ(\sqrt{d}_{1}, \sqrt{d}_{2})$, where $d_1$ and $d_2$ are  two coprime square-free integers such that  $d_1$ or $d_2\equiv1\pmod4$. Let $d_3=d_1d_2$ and assume $N(\epsilon_{3})=-1$,  $N(\epsilon_{1})= N(\epsilon_{2})=1$ or $N(\epsilon_{1})\neq N(\epsilon_{2})$. Let $i\neq j=1, 2$. Then
		$K$ is a P\'{o}lya field if and only if  one of the following conditions holds.
		\begin{enumerate}[\rm1.]
			\item\label{1i} $d_i=p_1$,  $d_j=p_2p_3$,
			\item\label{1j} $d_i=2$, $d_j=p_1p_2$,
			\item\label{1k} $d_i=2p_1$, $d_j=p_2$,
			\item\label{1l} $d_1=p_1p_2$,  $d_2=p_3p_4$,
			\item\label{1m} $d_i=2p_1$,  $d_j=p_2p_3$.
		\end{enumerate}
	\end{them}
	
	\begin{proof}
		As $N(\epsilon_{3})=-1$, so $d_3$ is not divisible by  any prime congruent to $3\pmod4$. Hence the  prime divisors of the discriminant $d_k$ of $k$ are exactly the prime divisors of  $d_3$. According to Theorem \ref{11}, we discuss two cases, let $i\neq j=1, 2$.
		\begin{enumerate}[\rm1.]
			\item Assume $N(\epsilon_{i})\neq N(\epsilon_{j})$ and $N(\epsilon_{3})=-1$. So by Theorem \ref{11} we have $\mathcal{P}_O(K)  \simeq E_{t-3}$. Then $K$ is a  P\'{o}lya field if and only if $t=3$. Thus $d_3=p_1p_2p_3$ or $2p_1p_2$.
			Suppose $(d_{i}, d_{j})\equiv(1, 1)\pmod4$,  then $d_i=p_1$ and $d_j=p_2p_3$. As in this case $N(\epsilon_{i})=-1$, so we must have $N(\epsilon_{j})=1$, on the other hand   $N(\epsilon_{3})=-1$.
			Suppose $(d_{i}, d_{j})\equiv (2, 1) \pmod4$,  then   $(d_i=2p_1$ and  $d_j=p_2)$ or $(d_i=2$ and  $d_j=p_1p_2)$. As $N(\epsilon_{i})\neq N(\epsilon_{j})$ and $N(\epsilon_{3})=-1$.
			\item Assume finally that $N(\epsilon_{1})=N(\epsilon_{2})=1$ and $N(\epsilon_{3})=-1$. So according to Theorem \ref{11} we have $\mathcal{P}_O(K)  \simeq E_{t-4}$. Then $K$ is a  P\'{o}lya field if and only if $t=4$. Thus $d_3=p_1p_2p_3p_4$ or $2p_1p_2p_3$.
			Suppose $(d_{i}, d_{j})\equiv(1, 1)\pmod4$,  thus  $d_i=p_1p_2$ and $d_j=p_3p_4$. Now if $(d_{i}, d_{j})\equiv (2, 1) \pmod4$,  then   $d_i=2p_1$ and $d_j=p_2p_3$.
		\end{enumerate}
		
	\end{proof}

	Now we apply Theorem \ref{19} to get the following theorems. Keep the notations adopted in the beginning of this section. To avoid long theorems, we will consider as above two cases, in the first one we assume  $N(\epsilon_{1})\neq N(\epsilon_{2})$ and in the second  we assume  $N(\epsilon_{1})= N(\epsilon_{2})=1$.
	\begin{them}
		Let $K=\QQ(\sqrt{d}_{1}, \sqrt{d}_{2})$,  where $d_1$ and $d_2$ are two coprime square-free integers such that $d_1$ or $d_2\equiv1\pmod4$. Let $k=\QQ(\sqrt{d_3})$, where $d_3=d_1d_2$, and assume that $d_3$ is divisible by a prime  $q\equiv 3\pmod4$ and $N(\epsilon_{1})\neq N(\epsilon_{2})$.  Let $i\neq j= 1, 2$ and denote  by $\epsilon_i$ the fundamental unit of $k_i$ such that $N(\epsilon_{i})=1$.
		\begin{enumerate}[\rm1.]
			\item\label{20} Assume  $\sqrt{\epsilon_i^a\epsilon_3}\in K$, for some $a=0, 1$, then $K$ is a P\'{o}lya field  if and only if   one of the following assertions holds.
			\begin{enumerate}[\rm i.]
				\item $d_i=q_1q_2$ and $d_j=p$,
				\item  $d_i=q_1q_2$ and $d_j=2$,
				\item $d_i=q$ and $d_j=p$,
				\item  $d_i=2q$ and $d_j=p$,
			\end{enumerate}
			\item  Assume  $\sqrt{\epsilon^a_i\epsilon_3}\not\in K$, for any $a=0, 1$,  then $K$ is a P\'{o}lya field  if and only if   one of the following assertions holds.
			\begin{enumerate}[\rm i.]
				\item  $d_i=p_1q$ and $d_j=p_2$,
				\item  $d_i=p_1q_1q_2$ and $d_j=p_2$,
				\item  $(d_i=q$ or $q_1q_2$ or $2q)$, $d_j=p_1p_2$,
				\item $d_i=pq_1q_2$ and $d_j=2$,
				\item $d_i=2q_1q_2$ and $d_j=p$,
				\item $d_i=2p_1q$ and $d_j=p_2$,
				\item $d_i=q_1q_2$ and $d_j=2p$.
			\end{enumerate}
		\end{enumerate}
	\end{them}
	
	\begin{proof}
		Assume $N(\epsilon_{1})\neq N(\epsilon_{2})$. So according to Theorem \ref{19} we have to consider  two cases. Recall that $d_3$ must be divisible by a prime number congruent to $3\pmod4$, so in this case, it may be that the number of the prime divisors of $d_k$ is not the same for $d_3$.
		\begin{enumerate}[\rm1.]
			\item If $\sqrt{\epsilon_i^a\epsilon_3}\in K$, for some $a=0, 1$,  then by Theorem  \ref{19} we have $\mathcal{P}_O(K)\simeq E_{t-3}$,  where $t$ is the number of prime divisors of $d_k$ the discriminant of $k$.  So $K$ is a  P\'{o}lya field if and only if $t=3$. Then we have the following  cases:\\
			If $d_3\equiv1\pmod4$, then $d_3=pq_1q_2$, i.e., $(d_i=q_1q_2$ and $d_j=p)$.\\
			If $d_3\equiv3\pmod4$, then $d_3=pq$ (remark that $t=3$). Thus $d_j=p$ and $d_i=q$.\\
			If $d_3\equiv2\pmod4$, then $d_3=2pq$ or $2q_1q_2$. Thus  $(d_i=2q$ and $d_j=p)$, $(d_i=2$ and $d_j=pq)$, $(d_i=q_1q_2$ and $d_j=2)$ or  $(d_i=2p$ and $d_j=q)$. For the last case to have $N(\epsilon_{i})=-1$.
			\item If $\sqrt{\epsilon_i^a\epsilon_3}\not\in K$, for any $a=0, 1$, then by Theorem  \ref{19} again  we have $\mathcal{P}_O(K)\simeq E_{t-4}$. So $K$ is a  P\'{o}lya field if and only if $t=4$. Taking into account that $N(\epsilon_{1})\neq N(\epsilon_{2})$, $(d_1\text{ or }d_2\equiv1\pmod4)$ and $d_3$ is divisible by a prime congruent to $3\pmod4$, and distinguishing  the cases $d_3\equiv1\pmod4$,  $d_3\equiv3\pmod4$ and  $d_3\equiv2\pmod4$,  we deduce all the possible cases.
		\end{enumerate}
	\end{proof}
	Now we treat the second case.
	\begin{them}
		Let $K=\QQ(\sqrt{d}_{1}, \sqrt{d}_{2})$,  $d_1$ and $d_2$ be two square-free integers such that $d_1$ or $d_2\equiv1\pmod4$ and $(d_{1}, d_{2})=1$.
		Let $k=\QQ(\sqrt{d_3})$, where $d_3=d_1d_2$, and assume that $d_3$ is divisible by a prime  $q\equiv 3\pmod4$ and $N(\epsilon_{1})= N(\epsilon_{2})=1$. Let $i\neq j= 1, 2$.
		\begin{enumerate}[\rm1.]
			\item Assume $\sqrt{\epsilon^a_1\epsilon^b_2\epsilon_3}\in K$ for some   $a, b\in \{0, 1 \}$. Then $K$ is a P\'{o}lya field if and only if  one of the following assertions holds.
			\begin{enumerate}[\rm i.]
				\item $d_i=q_1q_2$,  $d_j=2q_3$ or $q_3$
				\item $d_i=q_1q_2$,  $d_j=q_3q_4$.
				\item $d_i=p_1p_2$,  $(d_j=q$ or $q_1q_2$ or $2q)$,
				\item $d_i=q_1q_2$,  $d_j=2p$.
			\end{enumerate}
			\item Assume $\sqrt{\epsilon^a_1\epsilon^b_2\epsilon_3}\not\in K$ for any   $a, b\in \{0, 1 \}$. Then $K$ is a P\'{o}lya field if and only if  one of the following assertions holds.
			\begin{enumerate}[\rm i.]
				\item $d_i=q_1q_2$ and  $(d_j=2pq$ or $2q_3q_4$ or $pq_3q_4)$.
				\item $d_i=pq_1q_2$ and  $d_j=2q$.
				\item $d_i=pq_1$ and  $d_j=q_2q_3$.
				\item $d_i=pq_1q_2$ and  $d_j=q_3$.
				\item $d_i=p_1p_2$,  $(d_j=pq$ or $pq_1q_2$ or  $2pq$ or $2q_1q_2)$,
				\item $d_i=p_1q_1q_2$,  $d_j=2p$,
				\item $d_i=p_1p_2p_3$,  $(d_j=q$ or $q_1q_2$ or $2q)$,
				\item $d_i=q_1q_2$ and  $d_j=2p_1p_2$.
				
			\end{enumerate}
		\end{enumerate}
	\end{them}

	\begin{proof}
		Since $N(\epsilon_{1})= N(\epsilon_{2})=N(\epsilon_{3})=1$, then according to Theorem \ref{19} we have two cases to distinguish. Recall that $d_3$ must be divisible by a prime number congruent to $3\pmod4$, so in this case, it may be that the number of the prime divisors of $d_k$ is not the same for $d_3$, recall also that $d_1$ or $d_2\equiv1\pmod4$.
		\begin{enumerate}[\rm1.]
			\item Assume $\sqrt{\epsilon^a_1\epsilon^b_2\epsilon_3}\in K$,  where  $a, b\in \{0, 1 \}$, then   Theorem  \ref{19} implies that $\mathcal{P}_O(K)\simeq E_{t-4}$,  with $t$ is the number of prime divisors of $d_k$.  So $K$ is a  P\'{o}lya field if and only if $t=4$.
			Taking into account that $N(\epsilon_{1})= N(\epsilon_{2})=1$, $(d_1\text{ or }d_2\equiv1\pmod4)$ and $d_3$ is divisible by a prime congruent to $3\pmod4$, and distinguishing  the cases $d_3\equiv1\pmod4$,  $d_3\equiv3\pmod4$ and  $d_3\equiv2\pmod4$,  we deduce all the possible cases.
			\item Assume $\sqrt{\epsilon^a_1\epsilon^b_2\epsilon_3}\in K$,  for any  $a, b\in \{0, 1 \}$, then Theorem  \ref{19} implies that $\mathcal{P}_O(K)\simeq E_{t-5}$. So $K$ is a  P\'{o}lya field if and only if $t=5$.
			Taking into account that $N(\epsilon_{1})= N(\epsilon_{2})=1$, $(d_1\text{ or }d_2\equiv1\pmod4)$ and $d_3$ is divisible by a prime congruent to $3\pmod4$, and distinguishing  the cases $d_3\equiv1\pmod4$,  $d_3\equiv3\pmod4$ and  $d_3\equiv2\pmod4$,  we deduce all the possible cases.
		\end{enumerate}
		
	\end{proof}

	\section{\bf The first cohomology group of units of some real biquadratic fields}
	
Let $K=\QQ(\sqrt{d_1}, \sqrt{d_2})$ be a real biquadratic number field such that  $(d_1, d_2)=1$ and  $d_1$ or $d_2\equiv1\pmod4$, assume   $K$ is an unramified quadratic extension of $k=\QQ(\sqrt{d_1d_2})$. Let  $G_K$ be the Galois group $\mathrm{Gal}(K/\QQ)$ and $E_K$ be the group of units of $K$. In this section,  we will show that the first cohomology group $H^1(G_K, E_K)$ of units of $K$ is isomorphic to $E_b$ with $b\in\{2, 3, 4, 5\}$. This result is consistent with Zantema's corollary 5.6 (\cite{Za}), which gives that the order of $H^1(G_K, E_K)$ is bounded by $2^5$. The basis of our work in this section is the Proposition \ref{zatm6}, it is this which will give us the structure of $H^1(G_K, E_K)$. Denote by $\epsilon_1$, $\epsilon_2$ and $\epsilon_3$ the fundamental units of $k_1=\QQ(\sqrt{d_1})$, $k_2=\QQ(\sqrt{d_2})$ and $k=\QQ(\sqrt{d_1d_2})$ respectively.

	\begin{them}
		Let $K=\QQ(\sqrt{d}_{1}, \sqrt{d}_{2})$, where $d_1$ and $d_2$ are two coprime square-free integers. Assume $N(\epsilon_{3})=-1$,
		then  $$H^1(G_K, E_K)  \simeq \left\{ \begin{array}{ll}
			E_{2}&\text{ if } N(\epsilon_{1})=N(\epsilon_{2})=-1\text{ and }  \sqrt{\epsilon_1\epsilon_2\epsilon_3}\in K, \\
			E_{3}& \text{ if } (N(\epsilon_{1})=N(\epsilon_{2})=-1\text{ and }  \sqrt{\epsilon_1\epsilon_2\epsilon_3}\notin K)\text{ or }N(\epsilon_{1})\neq N(\epsilon_{2}), \\
			E_{4}& \text{ if }  N(\epsilon_{1})=N(\epsilon_{2})=1.
		\end{array}
		\right.$$
		
	\end{them}
	\begin{proof}
		Let $k=\QQ(\sqrt{d_3})$, where $d_3=d_1d_2$. Since $N(\epsilon_{3})=-1$ and $(d_{1}, d_{2})=1$, so $d_1$ or $d_2\equiv1\pmod4$, thus  the prime divisors of $d_k$, the discriminant of $k$, are exactly  those of $d_K$, the  discriminant of $K$. Note   also that $e_2$,  the ramification index of  $2$ in $K$, equals to $1$ or $2$.  As  $K/\QQ$ is a Galois extension,  so according to Proposition \ref{zatm6},  we have
		$$|H^1(G_K, E_K)\times||\mathcal{P}_O(K)|=\prod_{p|d_K}e_p,$$
		where $e_p$ is the ramification index of  $p$ in $K$. Thus
		$|H^1(G_K, E_K)|=2^t|\mathcal{P}_O(K)|^{-1}$, where $t$ is the number of prime divisors of $d_K$. By Theorem \ref{11} we have the following cases:
		\begin{enumerate}[\rm1.]
			\item If $N(\epsilon_{1})=N(\epsilon_{2})=-1$ and  $\sqrt{\epsilon_1\epsilon_2\epsilon_3}\in K$, then  $\mathcal{P}_O(K)\simeq E_{t-2}$. Thus
			$|H^1(G_K, E_K)|=2^t\times 2^{2-t}=2^2$. From which  we deduce that
			$H^1(G_K, E_K)\simeq E_2$.
			\item If $N(\epsilon_{1})=N(\epsilon_{2})=-1$ and  $\sqrt{\epsilon_1\epsilon_2\epsilon_3}\notin K$, then $\mathcal{P}_O(K)\simeq E_{t-3}$. Thus $|H^1(G_K, E_K)|=2^t\times 2^{3-t}=2^3$. From which  we deduce that  $H^1(G_K, E_K)\simeq E_3$.
			\item If $N(\epsilon_{1})\neq N(\epsilon_{2})$,  then $\mathcal{P}_O(K)\simeq E_{t-3}$. Thus $H^1(G_K, E_K)\simeq E_3$.
			\item Similarly, if $N(\epsilon_{1})=N(\epsilon_{2})=1$, then $\mathcal{P}_O(K)\simeq E_{t-4}$. Thus $H^1(G_K, E_K)\simeq E_4$.
			
		\end{enumerate}		
	\end{proof}
	We proceed similarly to prove the following results. Basing on Theorem \ref{19}, we obtain the following theorem.
	\begin{them}
		Let $K=\QQ(\sqrt{d}_{1}, \sqrt{d}_{2})$, where $d_1$ and $d_2$ are two coprime square-free integers
		satisfying $d_1$ or $d_2\equiv1 \pmod 4$. Let $k = \QQ(\sqrt{d_3})$, with $d_3 = d_1d_2$, and assume that $d_3$ is
		divisible by a prime $ \equiv3 \pmod 4$.
		Then  $$ H^1(G_K, 	E_K)  \simeq
		\left\{ \begin{array}{ll}
			E_{3}&\text{if } N(\epsilon_{1})\neq N(\epsilon_{2})\text{ and }\sqrt{\epsilon_i^a\epsilon_3}\in K, \text{ where }a\in\{0, 1\},
			i=1, 2 \text{ and }N(\epsilon_{i})=1,  \\
			E_{4}& \text{if } (N(\epsilon_{1})=N(\epsilon_{2})=1\text{ and }\sqrt{\epsilon^a_1\epsilon^b_2\epsilon_3}\in K \text{ for some } a,  b \in \{0, 1\})\text{ or}\\
			&(N(\epsilon_{1})\neq N(\epsilon_{2}) \text{ and }\sqrt{\epsilon_i^a\epsilon_3}\notin K, \text{ for all }a\in\{0, 1\},
			i=1, 2 \text{ and } N(\epsilon_{i})=1),    \\
			E_{5}& \text{if } \ N(\epsilon_{1})=N(\epsilon_{2})=1\text{ and }\sqrt{\epsilon^a_1\epsilon^b_2\epsilon_3}\notin K, \text{ for  all } a,  b \in \{0, 1\}.
		\end{array}
		\right.$$
	\end{them}

\end{document}